\documentclass[11pt]{amsart}
\usepackage{amsmath,amsthm,epsfig,amssymb}
\usepackage{xypic}
\input xy
\title{\it On the Erd$\ddot{O}$s-Straus conjecture}
\author{Eugen J. Ionascu and Andrew Wilson}
\curraddr{Department of Mathematics\\ Columbus State University\\4225 University Avenue\\
Columbus, GA 31907}
\date{January $7^{th}$, 2010}
\begin{document}
\def\sms{\small\scshape}
\baselineskip18pt
\newtheorem{theorem}{\hspace{\parindent}
T{\scriptsize HEOREM}}[section]
\newtheorem{proposition}[theorem]
{\hspace{\parindent }P{\scriptsize ROPOSITION}}
\newtheorem{corollary}[theorem]
{\hspace{\parindent }C{\scriptsize OROLLARY}}
\newtheorem{lemma}[theorem]
{\hspace{\parindent }L{\scriptsize EMMA}}
\newtheorem{definition}[theorem]
{\hspace{\parindent }D{\scriptsize EFINITION}}
\newtheorem{problem}[theorem]
{\hspace{\parindent }P{\scriptsize ROBLEM}}
\newtheorem{conjecture}[theorem]
{\hspace{\parindent }C{\scriptsize ONJECTURE}}
\newtheorem{example}[theorem]
{\hspace{\parindent }E{\scriptsize XAMPLE}}
\newtheorem{remark}[theorem]
{\hspace{\parindent }R{\scriptsize EMARK}}
\renewcommand{\thetheorem}{\arabic{section}.\arabic{theorem}}
\renewcommand{\theenumi}{(\roman{enumi})}
\renewcommand{\labelenumi}{\theenumi}
\newcommand{\Q}{{\mathbb Q}}
\newcommand{\Z}{{\mathbb Z}}
\newcommand{\N}{{\mathbb N}}
\newcommand{\C}{{\mathbb C}}
\newcommand{\R}{{\mathbb R}}
\newcommand{\F}{{\mathbb F}}
\newcommand{\K}{{\mathbb K}}
\newcommand{\D}{{\mathbb D}}
\def\phi{\varphi}
\def\ra{\rightarrow}
\def\sd{\bigtriangledown}
\def\ac{\mathaccent94}
\def\wi{\sim}
\def\wt{\widetilde}
\def\bb#1{{\Bbb#1}}
\def\bs{\backslash}
\def\cal{\mathcal}
\def\ca#1{{\cal#1}}
\def\Bbb#1{\bf#1}
\def\blacksquare{{\ \vrule height7pt width7pt depth0pt}}
\def\bsq{\blacksquare}
\def\proof{\hspace{\parindent}{P{\scriptsize ROOF}}}
\def\pofthe{P{\scriptsize ROOF OF}
T{\scriptsize HEOREM}\  }
\def\pofle{\hspace{\parindent}P{\scriptsize ROOF OF}
L{\scriptsize EMMA}\  }
\def\pofcor{\hspace{\parindent}P{\scriptsize ROOF OF}
C{\scriptsize ROLLARY}\  }
\def\pofpro{\hspace{\parindent}P{\scriptsize ROOF OF}
P{\scriptsize ROPOSITION}\  }
\def\n{\noindent}
\def\wh{\widehat}
\def\eproof{$\hfill\bsq$\par}
\def\ds{\displaystyle}
\def\du{\overset{\text {\bf .}}{\cup}}
\def\Du{\overset{\text {\bf .}}{\bigcup}}
\def\b{$\blacklozenge$}

\def\eqtr{{\cal E}{\cal T}(\Z) }
\def\eproofi{\bsq}

\begin{abstract} {\it Paul Erd$\ddot{o}$s conjectured that
 for every $n\in \mathbb N$, $n\ge 2$, there exist $a$, $b$, $c$
natural numbers, not necessarily distinct, so that
$\frac{4}{n}=\frac{1}{a}+\frac{1}{b}+\frac{1}{c}$ (see \cite{rg}).
In this paper we prove an extension of Mordell's theorem and
formulate a conjecture which is stronger than Erd$\ddot{o}$s'
conjecture.}
\end{abstract} \maketitle
\section{INTRODUCTION}
The subject of Egyptian fractions (fractions with numerator equal to
one and a positive integer as its denominator) has incited the minds
of many people going back for more than three millennia and
continues to interest mathematicians to this day. For instance, the
table of decompositions of fractions $\frac{2}{2k+1}$ as a sum of
two, three, or four unit fractions found in the Rhind papyrus has
been the matter of wander and stirred controversy for some time
between the historians. Recently, in \cite{aaa}, the author proposes
a definite answer and a full explanation of the way the
decompositions were produced. Our interest in this subjected started
with finding decompositions with only a few unit fractions.

 It is known that every positive rational number can be written as a
finite sum of different unit fractions. One can verify this using
the so called Fibonacci method and the formula
$\frac{1}{n}=\frac{1}{n+1}+\frac{1}{n(n+1)}$, $n\in \mathbb N$. For
more than three fourths of the natural numbers $n$, $\frac{4}{n}$
can be written as sum of only two unit fractions: the even numbers,
and the odd numbers $n$ of the form $n=4k-1$, via the identities
$\frac{2}{2k-1}=\frac{1}{k}+\frac{1}{k(2k-1)}$, and
$\frac{4}{4k-1}=\frac{1}{k}+\frac{1}{k(4k-1)}$, $k\in \mathbb N$. It
is clear that if we want to write

$$\frac{4}{n}=\frac{1}{a}+\frac{1}{b}+\frac{1}{c},\ \ a,b,c\in\mathbb N,$$

\n we can just look at primes $n$. However, just as a curiosity, for
$n=2009=7^2(41)$ (multiple of four plus one) one needs still two
unit fractions, and there are only three such representations:

$$\frac{4}{2009}=\frac{1}{504}+ \frac{1}{144648}=
\frac{1}{574}+\frac{1}{4018}=\frac{1}{588} +\frac{1}{3444}.$$

This follows from the following characterization theorem which is
well known (see \cite{ec} and \cite{gm}).

\begin{theorem}\label{twofractions} Let $m$ and $n$ two coprime
positive integers. Then $$\frac{m}{n}=\frac{1}{a}+\frac{1}{b}$$ for
some positive integers $a$ and $b$, if and only if there exists
positive integers $x$ and $y$ such that

(i) $xy$ divides $n$, and

(ii) $x+y\equiv 0$ (mod m).
\end{theorem}

\n In what follows we will refer to the equality

\begin{equation}\label{esconj}
\frac{4}{n}=\frac{1}{a}+\frac{1}{b}+\frac{1}{c}.
\end{equation}

\n and say that {\it $n$ has a representation} as in (\ref{esconj}),
or that (\ref{esconj}) {\it has a solution}, if there exist $a$,
$b$, $c$ ($a\le b\le c$) natural numbers, not necessarily different
satisfying (\ref{esconj}). Since $\frac{4}{n}> \frac{1}{a}$, the
smallest possible value of $a$ is $\lceil n/4\rceil$. The biggest
possible value of $a$ is $\lfloor \frac{3n}{4}\rfloor$, for instance
$\frac{4}{9}=\frac{1}{6}+\frac{1}{6}+\frac{1}{9}$.

If $n=4k+1$ ($k\in\mathbb N$) then we can try to use the smallest
value first for $a$, i.e. $a=k+1$:

\begin{equation}\label{firstequation}
\frac{4}{4k+1}=\frac{1}{k+1}+\frac{3}{(k+1)(4k+1)}.
\end{equation}

Now, if the second term in the right hand side of
(\ref{firstequation}) could be written as a sum of two unit
fractions we would be done. This is not quite how the things are in
general, but if we analyze the cases $k=3l+r$ with $r\in \{0,1,2\}$,
$l\ge 0$, $l\in \mathbb Z$, we see that there is only one excepted
case in which we get stuck: $k=3l$. This because
Theorem~\ref{twofractions} can be used in one situation: $k=3l+1$
implies $1+(3l+1+1)\equiv 0$ (mod 3). On the other hand, if $k=3l+2$
we get $(k+1)=(3l+2)+1=3(l+1)\ge 3$ and the second term is already a
unit fraction.

In order to simplify the statements of some of the facts in what
follows we will introduce a notation. For every $i\in \mathbb N$ let
$\cal C_i$ be defined by

$${\cal C}_i:=\{n | \  (\ref{esconj})\ has\ a\ solution\ with\ a\le \frac{n+4i-1}{4}\}.$$

It is clear that ${\cal C}_i\subset {\cal C}_{i+1}$ and then {\it
Erd$\ddot{o}$s}-Straus' conjecture is equivalent to $\underset{i\in
\mathbb N}{\bigcup} {\cal C}_i=\mathbb N$. Thus we obtained a pretty
simple fact about the Diophantine equation (\ref{esconj}):

\begin{proposition} \label{fprop} The equation (\ref{esconj}) has at least one
solution for every prime number $n$, except possibly for those
primes of the form $n\equiv 1$ (mod 12). Moreover,

$$\mathbb N\setminus {\cal C}_1\subset \{n| n\equiv 1\ (mod\ 12)\}.$$
\end{proposition}

We observe that $12=2^2(3)$, a product of a combination of the first
two primes. The first prime that is excluded in this proposition is
$13$. The equality (\ref{firstequation}) becomes

\begin{equation}\label{firstequationprim}
\frac{4}{12l+1}=\frac{1}{3l+1}+\frac{3}{(3l+1)(12l+1)}.
\end{equation}

At this point we can do another analysis modulo any other number as
long we can reduce the number of possible situations for which we
cannot say anything about the decomposition as in (\ref{esconj}). It
is easy to see that $3l+1$ is even if $l$ is odd and then
Theorem~\ref{twofractions} can be used easily with $x=1$ and $y=2$.
This means that we have in fact an improvement of the
Proposition~\ref{fprop}:

\begin{proposition} \label{sprop} The equation (\ref{esconj}) has at least one
solution for every prime number $n$, except possibly for those
primes of the form $n\equiv 1$ (mod 24). In fact,

$$\mathbb N\setminus {\cal C}_1\subset \{n| n\equiv 1\ (mod\ 24)\}.$$
\end{proposition}

Let us observe that $24+1=5^2$, $48+1=7^2$, which pushes the first
prime excluded by this last result to $73$. Quite a bit of progress
if we think in terms of the primes in between that have been taken
care off, almost by miracle.

If $n=24k+1$, then the smallest possible value for $a$ is $6k+1$ and
at this point let us try now the possibility that
$a=6k+2=\frac{n+7}{4}$,

\begin{equation}\label{equation2}
\frac{4}{24k+1}=\frac{1}{6k+2}+\frac{7}{2(3k+1)(24k+1)},\ \ k\in
\mathbb N.
\end{equation}

In the right hand side of (\ref{equation2}), the second term has a
bigger numerator but the denominator has now at least three factors.
This increases the chances that Theorem~\ref{twofractions} can be
applied and turn that term into a sum of only two unit fractions.
Indeed, for  $k=7l+r$, we get that $n=24k+1\equiv 0$ (mod 7) if
$r=2$, $2(3k+1)+1\equiv 0$ (mod 7) if $r=3$, $n+1=2(12k+1)\equiv 0$
(mod 7) if $r=4$, and $n+2=24k+3\equiv 0$ (mod 7) if $r=6$.
Calculating the residues modulo 168 in the cases $r\in \{0,1,5\}$ we
obtain:

\begin{proposition}\label{aftermod7}
The equation (\ref{esconj}) has at least one solution for every
prime number $n$, except possibly for those primes of the form
$n\equiv r$ (mod 168), with $r\in \{ 1, 5^2, 11^2\}$, $k\in \mathbb
Z$, $k\ge 0$. More precisely,
$$\mathbb N\setminus {\cal C}_2\subset \{n| n\equiv 1, 5^2, 11^2\ (mod\ 168)\}.$$
\end{proposition}

Let us observe that $168=2^3(3)(7)$, $168+1=13^2$, and the excepted
residues modulo 168 are all perfect squares. Because of this,
somehow, the first prime that is excluded by this result is
$193=168+25$. Again, we have even a higher jump in the number of
primes that have been taken care of. As we did before there is an
advantage to continue using (\ref{equation2}) and do an analysis now
on $k$ modulo $5$.

For $k=5l+r$, we have $n\equiv 0$ (mod 5) if $r=1$, $3k+1\equiv 0$
(mod 5) if $r=3$, and $6k+1\equiv 0$ (mod 5) if $r=4$, which puts
$n\in {\cal C}_2$ again. Therefore, we have for $r\in \{0,2\}$ the
following excepted residues modulo 120.

\begin{proposition}\label{aftermod5}
The equation (\ref{esconj}) has at least one solution for every
prime number $n$, except possibly for those primes of the form
$n\equiv r$ (mod 120), with $r\in \{ 1, 7^2\}$, $k\in \mathbb Z$,
$k\ge 0$. More precisely,
$$\mathbb N\setminus {\cal C}_2\subset \{n| n\equiv 1, 7^2\ (mod\ 120)\}.$$
\end{proposition}

One can put these two propositions together and get Mordell's
Theorem.

\begin{theorem}\label{mordel}
{\bf (Mordell \cite{gm})} The equation (\ref{esconj}) has at least
one solution for every prime number $n$, except possibly for those
primes of the form $n=840k+r$, where $r\in \{1, 11^2, 13^2, 17^2,
19^2, \ 23^2\}$, $k\in \mathbb Z$, $k\ge 0$. Moreover, we have

$$\mathbb N\setminus {\cal C}_2\subset \{n| n\equiv 1,  11^2, 13^2, 17^2, 19^2, 23^2\ (mod\ 840)\}.$$
\end{theorem}

\proof.\  By Proposition~\ref{aftermod7},  $n=168k+1$ may be an
exception but if $k=5l+r$, with $r\in \{0,1,2,3,4\}$ we have
$n\equiv 1$ or $7^2$ (mod 120) only for $r\in\{0,1\}$. These two
cases are the exceptions for both propositions and they correspond
to $n\equiv 1$ or $13^2$ (mod 840). All other excepted cases are
obtained the same way.\eproof

Let us  observe that $840=2^3(3)(5)(7)$ and the residues modulo 840
are all perfect squares. Not only that but $840+1=29^2$,
$840+11^2=31^2$, and $1009=840+13^2$ is the first prime that is
excluded by this important theorem. While 193 is the 44-th prime
number, 1009 is the $169^{th}$ prime. It is natural to ask if a
result of this type can be obtained for an even bigger modulo. We
will introduce here the next natural step into this analysis, which
implies to allow $a$ be the next possible value, i.e.
$\frac{n+11}{4}$, and we will be  using the identities


\begin{equation}\label{equations120}
\frac{4}{120k+1}=\frac{1}{30k+3}+\frac{11}{3(10k+1)(120k+1)},\ \
k\in \mathbb N,
\end{equation}

\begin{equation}\label{equations120k49}
\frac{4}{120k+49}=\frac{1}{30k+15}+\frac{11}{3(5)(2k+1)(120k+49)},\
\ k\in \mathbb N.
\end{equation}

\section{The analysis modulo 11 }

According to Proposition~\ref{aftermod5} we may continue to look
only at the two cases modulo 120 and use only the two formulae
above. If we continue the analysis modulo 11 in these two cases  we
obtain the following theorem.

\begin{theorem}\label{1320}
 The equation (\ref{esconj}) has at least
one solution for every prime number $n$, except possibly for those
primes of the form $n=1320k+r$, where $$r\in \{1, 7^2, 13^2, 17^2,
19^2, \ 23^2, 29^2, 31^2, 7(103), 1201, 7(127), 23(47)\}:=E, \ k\in
\mathbb Z, \ k\ge 0.$$ Moreover, we have

$$\mathbb N\setminus {\cal C}_3\subset \{n| n \in E (mod\ 1320)\}.$$
\end{theorem}

\proof.\ If $n=120k+1$ and $k=11l+1$, we see that $n\equiv 0$ (mod
11) and so (\ref{equations120}) gives the desired decomposition as
in (\ref{esconj}) right away. If $k=11l+r$ and $r\in \{2,4,5 \}$ the
Theorem~\ref{twofractions} can be employed to split the second term
in (\ref{equations120}) as a sum of two unit fractions. For
instance, for $r=5$ we have $1+3(10(11l+5)+1)\equiv 0$ (mod 11), so
one can take $m=11$, $x=1$ and $y=30k+3$ in
Theorem~\ref{twofractions}. Hence we have seven exceptions in this
situation:

\begin{itemize}
\item $r=0$ corresponds to $n\equiv 1$ (mod 1320),
\item $r=3$ gives $n\equiv 19^2$ (mod 1320),
\item $r=6$ corresponds to $n \equiv 7(103)$ (mod 1320),
\item $r=7$ gives $n\equiv 29^2$ (mod 1320),
\item $r=8$ corresponds to $n \equiv 31^2$ (mod 1320),
\item $r=9$ gives $n\equiv 23(47)$ (mod 1320), and finally
\item $r=10$ corresponds to $n \equiv 1201$ (mod 1320).
\end{itemize}

If $n=120k+49$ and $k=11l+r$, then for $r=5$ we have $n\equiv 0$
(mod 11). If $r\in \{3,5,6,8,9,10\}$ we can use
Theorem~\ref{twofractions}. The exceptions then are:

\begin{itemize}
\item $r=0$ corresponds to $n\equiv 7^2$ (mod 1320),
\item $r=1$ gives $n\equiv 13^2$ (mod 1320),
\item $r=2$ corresponds to $n \equiv 17^2$ (mod 1320),
\item $r=4$ gives $n\equiv 23^2$ (mod 1320),
\item $r=7$ corresponds to $n \equiv 7(127)$ (mod 1320). \eproof
\end{itemize}

Putting Theorem~\ref{aftermod7} and Theorem~\ref{1320} together we
get the following 36 exceptions:

$$\begin{array}{cccccc}
1^2  & 13^2& 17^2 & 19^2 & 23^2 & 29^2 \\
31^2 & \underline{1201} & 37^2 & 41^2 & 43^2 & 13(157)\\
  47^2 &\underline{2521}  &  19(139) & \underline{2689} & 53^2 & \underline{3361}\\
  59^2& \underline{3529}& 61^2 & 29(149) &   67^2&  71^2\\
13(397) & 73^2 & \underline{5569} & 17(353) & 31(199) & 79^2 \\
83^2& \underline{7561} & \underline{7681} & 89^2 & \underline{8089}
& \underline{8761}
\end{array}$$

The residue 1201, the first prime in this list is not really an
exception because of the following identity:

\begin{equation}\label{1201}
\frac{4}{9240k+1201}=\frac{1}{2310k+308}+\frac{1}{5(9240k+1)(15k+2)}+\frac{1}{770(9240k+1)(15k+2)},
\end{equation}

\n which shows that $9240k+1201\in {\cal C}_8$ for all $k\in \mathbb
Z$, $k\ge 0$.

We checked for similar identities and found just another similar
identity for the exception $17(353)=6001$:

\begin{equation}\label{6001}
\begin{array}{l}
\ds \frac{4}{9240k+6001}= \frac{1}{2310k+1540}+
\\ \\
\ds
\frac{1}{385(9240k+6001)(2034k+1321)}+\frac{1}{22(3k+2)(2034k+1321)},
\end{array}
\end{equation}

\n which shows that $9240k+6001\in {\cal C}_{40}$ for all $k\in
\mathbb Z$, $k\ge 0$.

\begin{theorem}\label{mordellextension}
The equation (\ref{esconj}) has at least one solution for every
prime number $n$, except possibly for those primes of the form
$n\equiv r$ (mod 9240) where $r$ is one of the 34 entries in the
table:

\vspace{0.1in}

 \centerline{
\begin{tabular}{|c|c|c|c|c|c|}
  \hline
  $1^2$  & $13^2$& $17^2$ & $19^2$ & $23^2$ & $29^2$ \\
   \hline
$31^2$  & $37^2$ & $41^2$ & $43^2$ & $13(157)$&  $47^2$ \\
 \hline
$\underline{2521}$  &  $19(139)$ & $\underline{2689}$ & $53^2$ & $\underline{3361}$ &  $59^2$\\
 \hline
$\underline{3529}$& $61^2$ & $29(149)$ &   $67^2$&  $71^2$ & $13(397)$\\
 \hline
 $73^2$ & $\underline{5569}$& $31(199)$ & $79^2$ & $83^2$ & $\underline{7561}$\\
  \hline
 $\underline{7681}$ & $89^2$ & $\underline{8089}$ &
 $\underline{8761}$
 & & \\
  \hline
\end{tabular}}

\vspace{0.1in}

Moreover, if $n$ is not of the above form, it is in the class ${\cal
C}_3$, or in ${\cal C}_8$ if $n\equiv 1201$ (mod 9240), or in ${\cal
C}_{40}$ if $n\equiv 6001$ (mod 9240).
\end{theorem}

\proof. \ We look to see for which values of $r\in
\{0,1,2,3,4,5,6\}$ we have $1320(7l+r)+s\in \{1,5^2,11^2\}$ (mod
168) with $s\in E$ ($E$ as in Theorem~\ref{1320}). For instance, if
$s=19^2$ we get that $r$ must be in the set $\{0,1,3\}$ in order to
have $1320(7l+r)+19^2\in \{1,5^2,11^2\}$ (mod 168). These three
cases correspond to residues $19^2$, $41^2$ and $29(149)$ modulo
9240. Each residue in $E$ gives rise to three exceptions. We leave
the rest of this analysis to the reader.\eproof

\section{Numerical Computations and Comments }

We observe that the first ten of these residues in
Theorem~\ref{mordellextension} are all perfect squares. In fact, all
19 squares of primes less than 9240 and greater than $11^2$ are all
excepted residues.  There is something curious about the fact that
all the perfect squares possible are excepted.  This may be related
with the result obtained by Schinzel in \cite{s} who shows that
identities such as (\ref{1201}), (\ref{6001}) and others in this
note, cannot exist if the residue is a perfect square. The same
phenomenon is actually captured in Theorem 2 in \cite{y}. The good
news about Theorem~\ref{mordel}, Theorem~\ref{1320}, and
Theorem~\ref{mordellextension}, is that the first excepted residues
are all perfect squares or composite and moreover their number is
essentially increasing with the moduli.

With our analysis unfortunately, there are a few other composite and
9 prime residues that have to be excluded. The prime 2521 is only
the 369-th prime and it is the first prime that is excluded by this
theorem. However, a decomposition  with the smallest $a$ possible is
exhibited in the equality

$$\frac{4}{2521} = \frac{1}{636} + \frac{1}{70588} +
\frac{1}{5611746},$$

\n which puts $2521\in {\cal C}_6$. The other primes are in the
smallest class $\cal C$ as follows:

\vspace{0.1in}

\centerline{
\begin{tabular}{|c|c|c|c|c|c|}
  \hline
  ${\cal C}_1$  & ${\cal C}_2$ & ${\cal C}_3$ & ${\cal C}_4$ & ${\cal C}_5$&  ${\cal C}_6$  \\
  \hline \hline \\
  3361, 7681, 8089 & 3529, 5569 & 8761 & 2689& 7561 & 2101, 2521  \\
  \hline
\end{tabular}}

\vspace{0.1in}

Clearly, one can continue this type of analysis by adding more
primes to the modulo which is at this point 9240. It is natural to
just add the primes in order regardless if they are of the form
$4k+1$ or $4k+3$. We see that {\it Erd$\ddot{o}$s}' conjecture is
proved to be true if one can show that the smallest excluded residue
for a set of moduli that converges to infinity is not a prime. One
way to accomplish this is to actually show that the pattern
mentioned above continues, i.e. the number of excluded residues
which are perfect squares or composite is essentially growing as the
modulus increases. This is actually our conjecture that we talked
about in the abstract. Numerical evidence points out that for
residues $r$ which are primes, we have $9240s+r\in {\cal
C}_{k(s,r)}$ with $k(s,r)$ bounded as a function of $s$. For
example, $9240s+2521 \in {\cal C}_{12}$ for every $s=1..100000$ and
the distribution through the smaller classes is

\vspace{0.1in}

\centerline{
\begin{tabular}{|c|c|c|c|c|c|c|c|c|c|c|c|}
  \hline
  ${\cal C}_1$  & ${\cal C}_2$ & ${\cal C}_3$ & ${\cal C}_4$ & ${\cal C}_5$&  ${\cal C}_6$ &
  ${\cal C}_7$& ${\cal C}_8$& ${\cal C}_9$ & ${\cal C}_{10}$ & ${\cal C}_{11 }$ &  ${\cal C}_{12 }$ \\
  \hline \hline \\
   10852 & 6444 & 5332 & 811 & 612 & 277 & 63 & 82 & 6 & 7 & 0 & 5 \\
  \hline
  $44.3\%$& $26.31\%$ & $21.78\%$ & $3.3\%$ & $2.5\%$ & $1.13\%$ & $0.26\%$ & $0.34\%$ & $0.025\%$
  & $0.029\%$  &  $0\%$ & $0.021\%$  \\
  \hline
\end{tabular}}
 \vspace{0.1in}
Now, if we add 13 to the factors we would have an analysis modulo
120120. It turns out that 2521 is not an modulo 120120 exception
because of the identity

$$\frac{4}{120120k+2521}=\frac{1}{30030k+4004}+\frac{1}{1001(120120k+2521)(810k+17)}
\frac{1}{22(15k+2)(810k+17)},$$

\n which shows that $120120k+2521\in {\cal C}_{3374}$ for all $k\in
\mathbb Z$, $k\ge 0$.

\n We found similar identities for the residues 2689, 3529, 29(149),
5569, 31(199), 7561, and 7681 modulo 120120. This suggests that one
may actually be able to obtain Mordell type results for bigger
moduli, in the sense that the perfect squares residues appear
essentially in bigger numbers, by implementing a finer analysis that
involves higher classes than $\cal C_3$. It is natural to believe
that this might be true taken into account that Vaughan \cite{v}
showed that

$$\frac{1}{m}\# \{n\in \mathbb N|\  n\le m,\  and\ ( \ref{esconj})\ does \ not\ have \ a\  solution
\}\le  e^{-c(\ln\ m)^{2/3}}, m\in \mathbb N,$$

\n for some constant $c>0$. This is saying, roughly speaking,  that
the proportion of the those $n\le m$ for which a writing with three
unit fractions of $4/n$ goes to zero a little slower than
$\frac{1}{m}$ as $m\to \infty$. The first few primes that require a
bigger class than the ones before are 2, 73, 1129, 1201, 21169,
118801, 8803369,..., corresponding to classes ${\cal C}_1$, ${\cal
C}_2$, ${\cal C}_3$, ${\cal C}_4$, ${\cal C}_8$, ${\cal C}_{15}$,
${\cal C}_{27}$, .... which shows a steep increase in the size of
classes relative to the number of jumps.

In \cite{y}, Yamamoto has a different approach to ours and obtains a
lesser number of exceptions at least for the primes involved in
Theorem~\ref{mordellextension}. For each prime $p$ of the form
$4k+3$ between 11 and 97, there is a table in \cite{y} of exceptions
for congruency classes $r$ ($n\equiv r$ (mod p) )  that is used to
check the conjecture using a computer for al $n\le 10^7$. However,
in \cite{rg}, Richard Guy mentions that the conjecture is checked to
be true for all $n\le 1003162753$.

We extended the search for a counterexample further for all $n\le
4,146,894,049$. For our computations we wrote a program that pushes
the analysis for a modulus of
$M=2,762,760=2^3(3)(5)(7)(11)(13)(23)$. The primes chosen here are
optimal, in the sense that the excepted residues are in number less
than the ones obtained by other options.   The first 12 exceptions
in this case are $1$, $17^2$, $19^2$, $29^2$, $31^2$, $37^2$,
$41^2$, $43^2$, $47^2$, $53^2$, $3361$, and $59^2$. The number of
these exceptions was 2299 but it is possible that our program was
not optimal from this point of view. Nevertheless, this meant that
we had to check the conjecture, on average, for every other $\approx
1201$ integer. The primes generated, 889456 of them, are classified
according to the smallest class they belong to in the next tables:

\vspace{0.1in}

\centerline{
\begin{tabular}{|c|c|c|c|c|c|c|c|c|c|c|c|}
  \hline
  ${\cal C}_1$  & ${\cal C}_2$ & ${\cal C}_3$ & ${\cal C}_4$ & ${\cal C}_5$&  ${\cal C}_6$ &
  ${\cal C}_7$& ${\cal C}_8$& ${\cal C}_9$ & ${\cal C}_{10}$ & ${\cal C}_{11 }$ &  ${\cal C}_{12 }$ \\
  \hline \hline \\
   380547 & 228230 & 128494& 61129& 50853& 17116 & 8459& 9580& 1836& 1386 & 547& 855\\
  \hline
  $42.8\%$& $25.7\%$ & $14.4\%$ & $6.9\%$ & $5.7\%$ & $2\%$ & $0.9\%$ & $1\%$ & $0.2\%$
  & $0.15\%$  &  $0.06\%$ & $0.096\%$  \\
  \hline
\end{tabular}}

\vspace{0.1in}

\centerline{
\begin{tabular}{|c|c|c|c|c|c|c|c|c|c|c|c|}
  \hline
  ${\cal C}_{13}$  & ${\cal C}_{14}$ & ${\cal C}_{15}$ & ${\cal C}_{16}$ & ${\cal C}_{17}$&  ${\cal C}_{18}$ &
  ${\cal C}_{19}$& ${\cal C}_{20}$& ${\cal C}_{21}$ & ${\cal C}_{22}$ & ${\cal C}_{23 }$ &  ${\cal C}_{24 }$ \\
  \hline \hline \\
   115 & 124 & 111 & 26 & 10 & 27 & 2 & 4 & 4 & 0 & 0 & 0 \\
  \hline
  $0.013\%$& $0.014\%$ & $0.012\%$ & $0.003\%$ & $0.001\%$ & $0.003\%$ & $0.0002\%$ & $0.00045\%$ & $0.00045\%$
  & $0\%$  &  $0\%$ & $0\%$  \\
  \hline
\end{tabular}}
 \vspace{0.1in}
\centerline{
\begin{tabular}{|c|c|c|}
  \hline
  ${\cal C}_{25}$  & ${\cal C}_{26}$ & ${\cal C}_{27}$  \\
  \hline \hline \\
   0& 0 & 1 \\
  \hline
  $0\%$& $0\%$ & $0.0001\%$  \\
  \hline
\end{tabular}}
 \vspace{0.1in}
 \vspace{0.1in}
So far, we have not seen a prime in a class ${\cal C}_{k}$ with $k>
27$. The result obtained in \cite{s} seems to imply that the minimum
class index for each prime, assuming the conjecture is true, should
have a limit superior of infinity.

\end{document}